\documentclass[11pt,reqno]{amsart}
\usepackage{amssymb,latexsym}
\usepackage{graphicx}
\usepackage{fancyhdr,amssymb}
\usepackage{url}		

\theoremstyle{plain}
\numberwithin{equation}{section}

\begin{document}
\fancyhead{}
\renewcommand{\headrulewidth}{0pt}
\fancyfoot{}
\fancyfoot[LE,RO]{\medskip \thepage}

\setcounter{page}{1}

\title[Fibonacci Numbers and Trivalent graphs ]{Fibonacci Numbers and Trivalent graphs}
\author{Cheng Lien Lang}
\address{Department Applied of Mathematics\\
                I-Shou University\\
                Kaohsiung, Taiwan\\
                Republic of China}
\email{cllang@isu.edu.tw}
\thanks{}
\author{Mong Lung Lang}
\address{Singapore
669608, Republic of Singapore   }
\email{lang2to46@gmail.com}

\begin{abstract}
We study the Fibonacci and Lucas numbers and demonstrate how identities can be constructed  by
 investigating trivalent graphs (an infinite graph with all vertices of degree 3) and splitting fields.
 \end{abstract}

\maketitle
\vspace{-.8cm}

\section{Introduction}
This article is a continuation of [7]. The main purpose of the present article is to demonstrate how
 identities can be constructed. The second section is devoted to the study of a trivalent graph associated
  to the Fibonacci numbers. It turns out that many identities (old and new) becomes {\em visible} by
   applying our technique. The third  section studies the splitting field of the characteristic polynomial
   of certain recurrences and such study again explains the existence of certain identities.
    The fourth  section studies the identities of the subsequences of the Fibonacci sequence.
    To the best of our knowledge, our  way of  finding the following third-order identity (see Proposition 2.5)
     is not in the literature yet.

  \medskip
  \noindent  {\bf Proposition 2.5.}
   {\small
       $$
       F_n^3 =
  ( F_{n-1}F_{n+1}F_{n+2}
    +F_{n-2}F_{n-1}F_{n+1})/6
    + (F_{n-2}F_{n+1}F_{n+2}
   -F_{n-2}F_{n-1}F_{n+2})/2
  . \eqno(1.1)$$}
  \noindent It seems that the technique we presented in the article will enable us to find
 more identities (see section 2.5).

\section{Fibonacci Numbers and trivalent graphs}

\noindent  Let  $e_1, e_2, e_3 $ be  arbitrary vectors
 placed in the following
 trivalent graph and let   $x$ be the vector  given by
$$ x = 2(e_3 +e_2)-e_1.\eqno(2.1)$$ Such a vector  $x$ is said to be
 $\mathcal F$-generated by $e_3$, $e_2$, and $e_1$ (in this order).
  Applying (2.1),  we may construct an infinite sequence as follows :
$$e_1,e_2,e_3, e_4 = 2(e_3+e_2)-e_1, \cdots, e_{n+1} = 2(e_n +e_{n-1})-e_{n-1}, \cdots .    \eqno(2.2)$$

\bigskip
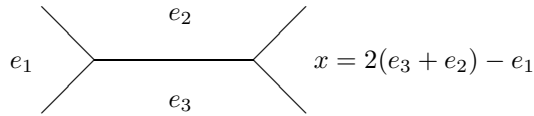
\begin{center}
 \begin{picture}(60,20)

 \multiput(2,0)(30,0){1}{\line(1,0){60}}

\multiput(2,0)(60,0){1}{\line(-1,1){20}}
 \multiput(2,0)(60,0){1}{\line(-1,-1){20}}

 \multiput(82,-20)(60,0){1}{\line(-1,1){20}}
 \multiput(82,20)(30,30){1}{\line(-1,-1){20}}

\put(-30, -3){\small  $e_1$}
\put(85, -3){\small $x = 2(e_3+e_2)-e_1$}
\put(30, 15){\small $ e_2$}
\put(30, -18){\small $e_3$}

 \end{picture}
\end{center}

\vspace{1cm}

\noindent
Denoted by $F(e_1,e_2,e_3)$ the above sequence.
In this section, we shall assign special  values for $e_1$, $e_2$ and $e_3$ and present a very
 visual method to display various identities concerning the Fibonacci numbers.
  Our work is inspired by Conway [2] that
        makes the values of binary quadratic forms visible. To put our idea
         in short, the trivalent graph we defined brings the Fibonacci numbers that
          reveal certain identities close to one another  and this makes such identities
           visible (see Proposition 2.2 for example).

\subsection{ Trivalent graph for $(1,0,0)$, $(2,1,0)$ and $(2,2,1)$}
 We shall first recall some results in [7].
A function $x(n)\, : \, \Bbb  Z\to \Bbb C$ is called an $\mathcal F$-function
 if $x(n)$ satisfies  the following recurrence relation.
 $$ x(n+3) = 2(x(n+2) + x(n+1)) - x(n).\eqno(2.3)$$
 One can show by induction that  $ x(n) = (-1)^n$, $x(n) = F_{n + r}^2$,
 $x(n) = L_{n + r}^2$,
   $x(n) = F_{n+r}F_{n+t}$,
     $x(n) = L_{n+r}L_{n+t}$,
    $x(n) = F_{n+r}L_{n+t}$,
    $x(n) = F_{2n+r}$  and  $x(n) = L_{2n+r}$
   are
    $\mathcal F$-functions,
    where $F_n$ and $L_n$ are  the $n$-th Fibonacci and Lucas numbers respectively.
      Further, if $A(n)$ and $B(n)$ are $\mathcal F$-functions, then
     $A(n) \pm B(n)$ are $\mathcal F$-functions. Note that composition and product
      of $\mathcal F$-functions are not $\mathcal F$-functions.
         The following is clear.

\medskip
\noindent {\bf Lemma 2.1.} {\em
 Let $A(n)$ and $B(n)$ be $\mathcal F$-function. Suppose that
  $A(k)= B(k)$ for
$k = 0,1$ and $2$. Then $A(n) = B(n)$.}

\medskip

\noindent {\bf Proposition 2.2.} {\em Let $F_n$ and  $L_n$ be the $n$-th Fibonacci and Lucas number. Then
$L_{n-1}^2 - F_{n-4}F_n- F_nF_{n+1} = F_{n-2}^2.$}

\medskip
\noindent {\em Proof.}  Recall first that $L_{-m} = (-1)^m L_m $
 and $F_{-m} = (-1)^{m+1}F_m$.
We shall assume  that $n \ge 0$. The case $n \le 0$ can be proved similarly.
Denoted by $A(n)$ and $B(n)$ the right and left hand side of our identity.
 Since both $A(n)$ and $B(n)$ are $\mathcal F$-functions  and $A(k)= B(k)$
  for $k= 0,1, 2$, we may apply Lemma 2.1  and conclude that
   $A(n) = B(n)$.\qed

\medskip
\noindent {\bf Discussion.} We shall now explain how to {\em see} Proposition 2,2 as follows.
 It is a trivalent graph that reveals the identity. Let $e_1=(1,0,0), e_2= (0,1,0), e_3 = (0,0,1)$.
The first nine terms of $F(e_1, 2e_1+e_2,2e_1+2e_2+e_3)$ are given as follows.

\bigskip
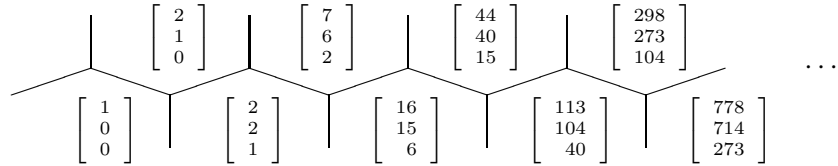
\begin{center}
 \begin{picture}(200,20)

\multiput(-30,0)(0,0){1}{\line(0,1){20}}
\multiput(-30,0)(0,0){1}{\line(-3,-1){30}}
\multiput(-30,0)(0,0){1}{\line(3,-1){30}}

\multiput(0,-10)(0,0){1}{\line(0,-1){20}}
\multiput(0,-10)(0,0){1}{\line(3,1){30}}

\multiput(30,0)(0,0){1}{\line(3,-1){30}}
\multiput(30,0)(0,0){1}{\line(0,1){20}}

\multiput(60,-10)(0,0){1}{\line(0,-1){20}}
\multiput(60,-10)(0,0){1}{\line(3,1){30}}


\multiput(90,0)(0,0){1}{\line(3,-1){30}}
\multiput(90,0)(0,0){1}{\line(0,1){20}}

\multiput(120,-10)(0,0){1}{\line(3,1){30}}
\multiput(120,-10)(0,0){1}{\line(0,-1){20}}

\multiput(150,0)(0,0){1}{\line(0,1){20}}
\multiput(150,0)(0,0){1}{\line(3,-1){30}}

\multiput(180,-10)(0,0){1}{\line(0,-1){20}}
\multiput(180,-10)(0,0){1}{\line(3,1){30}}



\put( -37, -25){\tiny$
\left [\begin{array} {r}
1\\0\\0\\
\end{array}\right ]
$}

\put( -9, 10){\tiny$
\left [\begin{array} {r}
2\\1\\0\\
\end{array}\right ]
$}

\put( 19, -25){\tiny$
\left [\begin{array} {r}
2\\2\\1\\
\end{array}\right ]
$}

\put( 47, 10){\tiny$
\left [\begin{array} {r}
7\\6\\2\\
\end{array}\right ]
$}

\put( 75, -25){\tiny$
\left [\begin{array} {r}
 16\\15\\6\\
\end{array}\right ]
$}

\put( 105, 10){\tiny$
\left [\begin{array} {r}
 44\\40\\15\\
\end{array}\right ]
$}

\put( 135, -25){\tiny$
\left [\begin{array} {r}
 113\\104\\40\\
\end{array}\right ]
$}

\put( 165, 10){\tiny$
\left [\begin{array} {r}
 298\\273\\104\\
\end{array}\right ]
$}

\put( 195, -25){\tiny$
\left [\begin{array} {r}
778\\714\\273\\
\end{array}\right ]
$}

\put (240, -2) {$\cdots$}

\end{picture}
\end{center}

\vspace{1.5cm}

\noindent One {\em sees} the following  interesting pattern from the first nine terms of our sequence.
\begin{enumerate}
\item [(i)] The second entry is a product of two Fibonacci numbers $F_nF_{n+1}$.
\item[(ii)] The difference of the first and second entries of every vector is the
 square of a Fibonacci number $F_{n-2}^2$.
 \end{enumerate}
(i) and (ii) of the above suggest us to express the first entry in terms of Fibonacci and /or Lucas numbers and we did.
it is given as follows.

\begin{enumerate}
\item [(iii)] The first entry of every vector is of the form $L_{n-1}^2 - F_{n-4}F_n$, where $L_k$
 is the $k$-th Lucas number.
\end{enumerate}

\noindent (i), (ii) and (iii) of the above implies that for the first nine vectors of the
 above graph, one has
$$L_{n-1}^2 - F_{n-4}F_n- F_nF_{n+1} = F_{n-2}^2,\eqno (2.4)$$
 which leads us the Proposition 2.2.
Note also that

\begin{enumerate}
\item[(iv)] The sum of the third entry (product of two Fibonacci numbers) of the first $2k+1$ terms is a square of a Fibonacci number.
    \item[(v)] the second and third entries of every vector are product of Fibonacci numbers and their sum
     is another Fibonacci number.
\end{enumerate}

\subsection
{Trivalent graph for $(1,2), (2,1)$ and $(1,2)$}
We shall present in this section two identities involved sum of Fibonacci numbers.
As the proof of such identities is elementary, we shall present how such identities  can be
    seen rather than how such identities can be proved.

\medskip
\noindent{\bf Proposition 2.3.}
\begin{enumerate}
\item[(i)]$3(F_4+ F_{12} + F_{20 }+ \cdots
 + F_{4+ 8n} ) = F_{4n+4}^2$.
 \item[(ii)]
 $  1+ 3(F_6+ F_{14} + F_{22 }+ \cdots
 + F_{6+ 8n} ) = F_{4n+ 5}^2.$
 \end{enumerate}

\medskip
\noindent {\bf Discussion.}
Let $e_1 = (1,0),$ $e_2=(0,1)$ and let $ u_1 =e_1+2e_2, u_2= 2e_1+e_2, u_3 = e_1+2e_2$
The first ten vectors $u_1,u_2,\cdots , u_{10}$ of $F(u_1,u_2,u_3)$
 are given as follows.   Note that
 $e_1 +2e_2$ appears twice in the beginning.

\bigskip
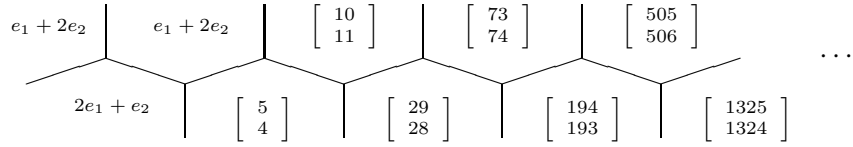
\begin{center}
 \begin{picture}(200,20)

\multiput(-30,0)(0,0){1}{\line(0,1){20}}
\multiput(-30,0)(0,0){1}{\line(-3,-1){30}}
\multiput(-30,0)(0,0){1}{\line(3,-1){30}}

\multiput(0,-10)(0,0){1}{\line(0,-1){20}}
\multiput(0,-10)(0,0){1}{\line(3,1){30}}

\multiput(30,0)(0,0){1}{\line(3,-1){30}}
\multiput(30,0)(0,0){1}{\line(0,1){20}}

\multiput(60,-10)(0,0){1}{\line(0,-1){20}}
\multiput(60,-10)(0,0){1}{\line(3,1){30}}


\multiput(90,0)(0,0){1}{\line(3,-1){30}}
\multiput(90,0)(0,0){1}{\line(0,1){20}}

\multiput(120,-10)(0,0){1}{\line(3,1){30}}
\multiput(120,-10)(0,0){1}{\line(0,-1){20}}

\multiput(150,0)(0,0){1}{\line(0,1){20}}
\multiput(150,0)(0,0){1}{\line(3,-1){30}}

\multiput(180,-10)(0,0){1}{\line(0,-1){20}}
\multiput(180,-10)(0,0){1}{\line(3,1){30}}

\put( -66, 10){\tiny$e_1+2e_2$}
\put( -12, 10){\tiny$e_1+2e_2$}
\put( -42, -20){\tiny$2e_1+e_2$}

\put( 18, -25){\tiny$
\left [\begin{array} {r}
5\\4\\
\end{array}\right ]
$}

\put( 47, 10){\tiny$
\left [\begin{array} {r}
10\\11\\
\end{array}\right ]
$}

\put( 75, -25){\tiny$
\left [\begin{array} {r}
29\\28\\
\end{array}\right ]
$}

\put( 105, 10){\tiny$
\left [\begin{array} {r}
73\\74\\
\end{array}\right ]
$}

\put( 135, -25){\tiny$
\left [\begin{array} {r}
 194\\193\\
\end{array}\right ]
$}

\put( 165, 10){\tiny$
\left [\begin{array} {r}
 505\\506\\
\end{array}\right ]
$}

\put( 195, -25){\tiny$
\left [\begin{array} {r}
1325\\1324\\
\end{array}\right ]
$}

\put (240, -2) {$\cdots$}

\end{picture}
\end{center}
\bigskip

\vspace {1cm}
\noindent One sees that

\begin{enumerate}
\item[(i)] The differences $d_i$'s of the entries of the  vectors
 $(u_5, u_3)$, $(u_9, u_7), $ $(u_{13}, u_{11}) , \cdots $ of the top half of the
graph are three times of some Fibonacci numbers. The first three differences are $10-1=9=3F_4$, $505-73=
432=3F_{12}$ and
$23761-3466=20295=3F_{20}$.
\item[(ii)] The sum of consecutive $d_i$'s are squares. $9=3^2$, $9+432 = 21^2$, $9+432
+20295 = 144^2.$
\end{enumerate}
\noindent One sees that (i) and (ii) give us
$$ 3F_4+3F_{12} = F_8^2, \,\,\,\,3F_4+3F_{12} +3F_{20} = F_{12}^2,\eqno(2.5)$$
which leads us to (i) of Proposition 2.3.

\medskip
\noindent We now observe  the bottom half  of the above graph. The entries of the first few vectors gives the following.
$$1+ 3F_6 +3F_{14} = (-1+2) + (-5+29)
 + (-194+1325)= 34^2 = F_9^2,\eqno(2.6)$$
  which leads us to (ii) of Proposition 2.3.

\subsection {Trivalent graph of $ (1,0,-1),(2,1,1)$ and $(2,1,4)$} Similar to sections 2.1 and 2.2,
 the following identity can be obtained by study the trivalent graph of
 $ (1,0,-1),(2,1,1)$ and $(2,1,4)$.

\medskip
\noindent {\bf Proposition 2.4.} {\em
$ L_{2n+1}-F_{n+1}^2- ( L_n^2-F_{n-3}F_{n+1}+F_{2n-2}) = (-1)^{n-1}.$}

\subsection{More identities.} We  list in this section  a few identities
 that can be {\em seen} by studying certain trivalent graphs.
Let  $u_1, u_2 , u_3$ be three vectors.  The
 identity generated by
 $F(u_1, u_2, u_3)$ is given as follows.

{\small
$\begin{array} {llrllrrrrr}

(c1)\,\,\,\,\,\,\,\, &(1,1,2)&(0,1,-1)\,\,\,& (1,0,2)\,\,\,\,\, & :& &&  F_{2n}= F_{n+1}^2-F_{n-1}^2\\

\\

(c2) &(1,1,2)&(0,1,-1)& (1,0,2) & :& &&   F_{2n+1}= F_{n+1}^2+F_{n}^2\\

\\

(c3) &(1,1,2)&(0,1,-1)& (1,0,2) & :& &&   F_{n+2}F_{n-1}= F_{n+1}^2 -F_{n}^2\\

\\

(c4) &(1,1,2)&(0,1,-1)& (1,0,2) & :& &&   F_{n}^2 - F_{n-2}F_{n+2} = (-1)^n\\

\\

(c5) &(0,1,0)&(1,0,1)& (1,2,3) & :& &&   \sum_{i=1}^n F_i^2 = F_n F_{n+1}\\

\\

(c6) &(0,1,0)&(1,0,1)& (1,2,3) & :& &&   F_n^2 -F_{n-1}F_{n+1} = (-1)^{n-1}\\

\\

(c7) &(1,-1,0)&(0,1,0)& (0,0,1) & :& &&  F_n F_{n+3} = F_{n+1} F_{n+2} +(-1)^{n-1}\\

\\

(c8)&(1,1,1)&(0,1,0)&(0,0,1) & :& &&  F_n^2 -F_{n-1}^2 = F_n F_{n-1} +(-1)^{n-1}\\
\\
(c9) &(-1,0,0)&(3,1,0)&(0,0,1)& :&\,\,\,\,\,\, &   \,\,\,\, &   F_{2n+1}+(-1)^n =F_{n-1}F_{n+1} +F_{n+1}^2\\
\\

\end{array}
$}

\medskip

\noindent {\em Proof.} Note first that the left and right hand side of (c1)-(c9) are all
 $\mathcal F$-functions.
 One may now prove the identities by applying Lemma 2.1.\qed

  \subsection {Discussion} We have presented in this  section how identities about Fibonacci
   numbers (old and new) can be seen by drawing a simple trivalent graph.
   Note that the graph we use to generate various identities here is just a subgraph of an
   infinite trivalent graph (see Appendix A of [7]). Such infinite trivalent  graph may generate
   more identities. We are currently investigating this infinite graph and the related
    identities. Note that
        the vectors $e_i$' in this article are vectors in $\Bbb R^2$ or $\Bbb R^3$
     and the nonidentity constant we use in the recurrence relation $e_n = c(e_{n-1}+e_{n-2})-e_{n-3}$
     is $c = 2$. All these invariance can be replaced by other choices.
      For instance,
      $$e_n = 5e_{n-1}+15e_{n-2}-15e_{n-3} - 5e_{n-4} + e_{n-5}.\eqno(2.7)$$
      Note that $F_n^2$ satisfies the recurrence (2.1) and $F_n^4$ satisfies the recurrence (2.7).
      In (2.7), let $e_1=(0,0,0, 1), e_2 = (0,0,1,0), e_3=(0,0,0,0)$,$e_4=(0,1,0,0)$
       and $e_5= (1,0,0,0)$. The first
       nine $e_i$'s are

\begin{center}
$\begin{array} {llrrrrrrrr}
\\
e_1& e_2&e_3&e_4 &e_5&e_6& e_7& e_8& e_9\\
\\
0 & 0 & 0& 0 & 1&5& 40  & 260 & 1820 \\
0 & 0 & 0& 1 & 0 &15& 60  & 520 & 3276 \\
0 & 1 & 0& 0 & 0 &-5& -24 & -195 & -1260 \\
1 &0 & 0 & 0& 0 & 1&5& 40  & 260  \\
\\
\end {array}   $
\end{center}

 \noindent
 One sees easily that for each $ e_i = (a,b, -c, d)$, $6a,$ $2 b$, $2c$ and $6d$ are product of four
  Fibonacci numbers. The expression $ a+b -c+d $ is a fourth power  of a Fibonacci number.
Similar to Propositions 2.2-2.4,
       one has the following third-order identity.

  \medskip
  \noindent
 \medskip
\noindent  {\bf Proposition 2.5.}
   {\small
       $$
       F_n^3 =
  ( F_{n-1}F_{n+1}F_{n+2}
    +F_{n-2}F_{n-1}F_{n+1})/6
    + (F_{n-2}F_{n+1}F_{n+2}
   -F_{n-2}F_{n-1}F_{n+2})/2
  . \eqno(2.8)$$}

\noindent  To the best of our knowledge,
the way we construct  (2.8)
 is not in the literature yet. Note that our
 approach towards the study of the identities is computer-free and  our technique always makes the identities
  visible. However,
  it is perhaps possible
   to implement  this idea in a computer to visualise more identities. The recurrence relation
    one should consider is (see [5])

     $$ e_{n+1} = -\sum_{r=1}^{n+1}(-1)^{r(r+1)/2}\frac{\,\,F_{n+1} F_{n}\cdots F_{n+2-r}\,\,}
     {F_rF_{r-1} \cdots F_1}e_{n+1-r}.\eqno(2.9)$$

\noindent The interested readers may want to apply the technique of induction and/or Binet's formula
to obtain
 a general formula for $F_n^m$.

\section {Characteristic polynomials and Splitting fields }

\noindent Let $\{u_0, u_1, u_2, \cdots\}$ be a sequence defined by the following recurrence relation
$u_0 = 0, u_1= 1$,
 $$u_r= pu_{r-1} -qu_{r-2}. \eqno(3.1)$$

 \medskip
 \noindent
 where $p$ and $q$ are rational numbers. The characteristic polynomial of $
  \{u_r^n\}$ is defined to be (see for examples, [1], [3], [4], [6])

 $$\Phi_n(p, q, x)=
  \sum_{i=0}^{n}
  (-1)^{i} q^{i(i-1)/2} (n|i)_u x^i,\eqno(3.2)$$

 \medskip
 \noindent
 where $(n|k)_u$ is the generalised binomial coefficient (see Appendix A). Note that
$(n|0)_u=1$ and $(n|k)_u= u_nu_{n-1}\cdots u_{n-k+1}/u_ku_{k-1}\cdots u_1$ for $1\le k \le n$
 if $u_1 u_2 \cdots u_k \ne 0$.
  Let $\sigma$ and $\tau$ be roots of $x^2- px+q =0$. It is well known that

   $$\Phi_{n}(p,q, x) = \prod_{j=0}^{n} (x-\sigma^j\tau^{n-j}).\eqno(3.3)$$

 \medskip
 \noindent
 In this section,
 we shall associated to the above factorisation some identities, which we propose to call them
  the {\em Galois Identities} of $\Phi_n(p,q,x)$.
These identities  (see Proposition 3.4 and Corollary 3.6) must be well known among the experts. However,
our interpretation of the existence of these identities maybe  of some interest (see Discussion 3.7).

 \subsection {Fibonacci and Lucas numbers and $L_n^2 - 5F_n^2 = 4(-1)^n$}

\noindent
In (3.1), $p=1$ and $q=-1$ give the Fibonacci numbers.
 Applying (3.3), it is clear that
  the Galois group of $\Phi_n(1, -1, x) $ over $\Bbb Q$ is $\Bbb Z_2$ and the
   splitting field of $\Phi_n(1,-1, x)$ is $\Bbb Q(\sqrt 5)$.
    Recall another well known fact about the factorisation of $\Phi_n(1,-1, x)$.

  $$\Phi_n(1,-1,x) = (-1)^n(x^2-L_nx +(-1)^n)\Phi_{n-2}(1,-1, -x).\eqno(3.4)$$

   \medskip
   \noindent
    Since $\Phi_n(1,-1,x)$ splits in $\Bbb Q(\sqrt 5)$ and $x^2-L_nx +(-1)^n =0$
   splits in $\Bbb Q(\sqrt{L_n^2 -4(-1)^n})$, one must have

   $$L_n^2 -4(-1)^n = 5A_n^2,\eqno(3.5)$$

\medskip
   \noindent for some $A_n\in \Bbb N$.
   This tells us that the difference between $4(-1)^n$ and the square of the Lucas number $L_n$ must be five times a square $A_n^2$. As for why $A_n$ must be $F_n$,
    we note that both $F_n^2$ and $L_n^2 $ satisfy the recurrence relation (2.3)
     and that     $L_n^2-4(-1)^n$ and $F_n^2$ have the same initial values (see Lemma 2.1).
      As a consequence,
 we have just recover the following well known identity by investigating the splitting
 field of $\Phi_n(1,-1,x)$.

 $$L_n^2 -4(-1)^n = 5F_n^2.\eqno(3.6)$$

 \medskip
 \noindent {\bf Discussion.} The above investigation suggests that (3.6)
 is not just a numerical coincidence and
 can  be viewed as the consequence of the fact that  the Galois group of $\Phi_n(1,-1,x)$
  is $\Bbb Z_2$ and that the splitting field of $\Phi_n(1,-1,x)$ is
  $\Bbb Q(\sqrt 5) = \Bbb Q(\sqrt { L_n^2 + 4(-1)^n})$.
\medskip

\noindent {\bf Remark.} It follows easily from (3.6) that the following identities hold,
 the first identity is proved by Freitag and the second by Zeitlin and Filipponi (independently).
  See  [4] for more detail.

  $$\frac{\,L_n^2 - (-1)^a L_{n+a}^2\,}{F_n - (-1)^a F_{n+a}^2} = 5,\,\,
  \frac{\,L_n^2 + L_{n+2a}^2+ 8(-1)^n\,}{F_n + F_{n+2a}^2} = 5.\eqno(3.7)$$

\subsection {The Galois Identities of $w_r = pw_{r-1}-qw_{r-2}$}
By (3.3),    the splitting field of $\Phi_n(p, q, x)$ is
     $\Bbb Q(\sqrt { p^2 - 4q})$
and
     the Galois group
     of $\Phi_n (p, q, x) $ is $  \Bbb Z_2$
     if and only if $p^2-4q$ is not a perfect square in $\Bbb Q$.
\medskip

\noindent {\bf Definition 3.1.} Let $\sigma$ and $\tau$ be given as
  in (3.3). Define $\{w_n\}$ to be the sequence  $w_n = \sigma^n +\tau^n$.

  \medskip
  \noindent {\bf Lemma 3.2.} {\em $w_0=2,$ $w_1= p$, $w_r= u_{r+1} -qu_{r-1}$ and
  $w_r= pw_{r-1}-qw_{r-2}$. Suppose that
   $p^2-4q\ne 0$. Then $u_r= (w_{r+1}-qw_{r-1})/(p^2-4q)$. }
   \medskip

   \noindent {\em Proof.} It is clear that $w_0=2$ and that $w_1=p$. Applying Binet's formula, one has
    $w_r= u_{r+1} -qu_{r-1}$ and that $w_r$  satisfies the recurrence
    $w_r= pw_{r-1} -qw_{r-2}$. Since

    $$w_r= u_{r+1} -qu_{r-1} = pu_r -2qu_{r-1}, \, w_{r-1} = u_r-qu_{r-2}= 2u_r -pu_{r-1},\eqno(3.8)$$

 \medskip
 \noindent we conclude that   $u_r= (w_{r+1}-qw_{r-1})/(p^2-4q)$. This completes the proof of the lemma. \qed

\medskip
\noindent {\bf Lemma 3.3.} {\em  $q^n$,  $w_n^2$ and $u_n^2$ satisfy the following recurrence

 $$X(m+3) = (p^2-q)X(m+2) +(q^2-p^2q)X(m+1) +q^3X(m).\eqno(3.9)$$}

\noindent
{\em Proof.} Applying (3.1), one can show easily that
 $u_n^2$ satisfies the recurrence (3.9). The rest can be verified
  similarly. \qed

\medskip
\noindent
Suppose that $p^2 - 4q$ is not a square in $\Bbb Q$.
 Applying Galois Theory, the set of conjugates of $\sigma^n$ over $\Bbb Q$ is $\{\sigma^n, \tau^n\}$.
 Hence the following holds for every $n \in \Bbb N$.
 $$ f_n(x) = (x-\sigma^n)(x-\tau^n)
 = x^2 - w_n x + q^n
 \in \Bbb Q[x].\eqno (3.10)$$

\medskip
\noindent
It is clear that $f_n(x)$ splits in $\Bbb Q(\sqrt { w_n^2-4q^n})$.
 Applying (3.3),  $f_n(x)$ splits in  $\Bbb Q(\sqrt {p^2-4q})$.
  Hence
  $$ w_n^2 - 4q^n = z^2(p^2-4q),\eqno(3.11)$$

  \medskip
  \noindent for some $z \in \Bbb Q$. The following proposition shows that the solution of (3.11) is $z = u_n$,
   the recurrence we defined in (3.1).

\medskip
\noindent {\bf Proposition 3.4.} {\em Let $ w_n $ be given as in Definition $3.1$. Then $ w_n^2 - 4q^n = u_n^2(p^2-4q).$}

\medskip
\noindent {\em Proof.}
Since both $ w_n^2 - 4q^n$   and $u_n^2(p^2-4q)$ satisfy the recurrence (3.9)
 and admit the same initial values, we have
$ w_n^2 - 4q^n = u_n^2(p^2-4q).$ \qed
\medskip

\noindent {\bf Corollary 3.5.} {\em
The equation $x^2 +y^2 -z^2 = 4$ is solvable in $\Bbb Z$.
Further, one may choose $y$ and $z$ in such a way that $py = 2z$ for any $p \in \Bbb Z$.
 }

\medskip

\noindent {\em Proof.} Let $q=1$ and let $p$ be any integer. Applying  Proposition 3.4, one has
$w_n^2 +(2u_n)^2 - (pu_n)^2 = 4$ for all $n \ge 1$.\qed

\medskip
\noindent {\bf Corollary  3.6.} {\em$ w_{2n}- 2q^n = u_n^2(p^2-4q).$
In particular, $L_{2n} -2(-1)^n = 5F_n^2$.
}

\medskip
\noindent {\em Proof.} Applying Definition 3.1, one has $w_n^2 = w_{2n} +2q^n$. \qed

\medskip
\noindent {\bf Discussion 3.7.}
The identities  in Proposition 3.4  and Corollary 3.6 must be  well known and we propose to call them the {\em Galois identities} associated to
 $\Phi_n(p, q, x)$. We would like to emphasise that  Proposition 3.4 is not just a numerical
  coincidence but can  be treated as the consequence of the fact
   that
the splitting field of $\Phi_n(p,q,x)$ is
  $\Bbb Q(\sqrt {p^2-4q}) = \Bbb Q(\sqrt { w_n^2 -4q^n})$. The most famous identity among all, of
   course, is (3.6).

\section{Subsequences of the Fibonacci sequence}

\noindent One sees easily that
$$F_{n+r} = L_r F_n +(-1)^{r+1} F_{n-r}.\eqno (4.1)$$

\medskip
\noindent
The above is known as the multiple angle recurrence (see [8]).
In the case $r=2$, (4.1) gives $F_{n+2} = 3 F_n - F_{n-2}.$ It follows that both
 $x_n = F_{2n}$ and $y_n= F_{2n+1}$ ($\{x_n\}$ and $\{y_n\}$ are subsequences of the Fibonacci sequence) satisfy the recurrence $ w_n= 3w_{n-1}- w_{n-2}$.
  Applying Lemma 3.3,

  $$w_n^2 = 8 w_{n-1}^2 -8 w_{n-2}^2 + w_{n-3}^2.\eqno (4.2)$$

 \noindent
 As a consequence,

  $$F_{2n +k}^2 = 8 F_{2(n-1)+k}^2 -8 F_{2(n-2)+k}^2 + F_{2(n-3)+k}^2
  ,\eqno (4.3)$$

\medskip
\noindent
 where $k = 0$ or 1. Note that  $r$ can be any integer and the idea
  can be generalised to higher order products.

\section {Appendix A}
  Let $\sigma$ and $ \tau$ be roots of $x^2 -px+q=0$ (known as the characteristic polynomial
  of (3.1)). It is well known that
 $$u_r = \sum_{i=1}^{r-1} \sigma^{r-1-i}\tau^i = f_r(\sigma, \tau) ,\eqno(A1)$$

 \medskip
\noindent where $f_r(x,y)$ is the polynomial  $ \sum x^{r-1-i}y^i  =(x^r-y^r)/(x-y).$
 Let $G(z)$ be the following function.
 $$ G(z) = (1-z^{m+n})(1-z^{m+n-1})\cdots (1-z^{m+1}/(1-z^n)(1-z^{n-1})\cdots (1-z).\eqno(A2)$$
 $G(z)$ is known as the Gaussian binomial coefficient  and is a polynomial in $z$ (a more powerful result actually implies that (A2) can be written
  as product of cyclotomic polynomials).
  One may apply this fact to show that
  that
 $(f_r f_{r-1}\cdots f_{r-k+1})/(f_kf_{k-1}\cdots f_1)\in \Bbb Z[x,y]$ is  a polynomial in $x$ and $y$.
Denoted by $F(r,k, x, y)$ this polynomial. Define

 $$ (r|k)_u = F(r,k,\sigma, \tau).\eqno(A3)$$

\medskip
\noindent We call $(A3)$ the generalised binomial coefficient. It is clear that if $u_1 u_2 \cdots u_k
 \ne 0$. Then $(n|k)_u$ takes the following form.

 $$(n|k)_u= u_nu_{n-1}\cdots u_{n-k+1}/u_ku_{k-1}\cdots u_1. \eqno(A4)$$

\bigskip

\noindent MSC2010 : 11B39, 11B83.

\medskip

\end{document}